\documentclass{article}

\usepackage{fancyhdr}
\newlength{\myparskip}              \newlength{\myparindent}             %
\newlength{\oldparskip}             \newlength{\oldparindent}            %
\parskip\myparskip                  \parindent\myparindent               %
 \fancypagestyle{plain}{%
     \fancyhead[LO]{}
     \fancyhead[CO]{}\fancyhead[RO]{}
 }

\usepackage{modamsthm,modnatbib}

\usepackage{latexsym,amssymb,lastpage}
\usepackage{graphicx,amsfonts}
\usepackage{times,mathptmx,bm,amsmath}
\usepackage{dcolumn}

\renewenvironment{proof}[1][Proof]{\noindent\textit{#1. } }{\hfill$\square$}

 \newtheoremstyle{theorem}{6pt}{6pt}{\rm}{}{\sffamily}{ }{ }{}
 \theoremstyle{theorem}

  \newtheoremstyle{thm}{6pt}{6pt}{\rm}{}{\sffamily}{ }{ }{}
 \theoremstyle{thm}
\newtheorem{thm}{\sc Theorem}[section]

 \newtheoremstyle{lemma}{6pt}{6pt}{\rm}{}{\sffamily}{ }{ }{}
 \theoremstyle{lemma}
 \newtheorem{lemma}{\sc Lemma}[section]
 \newtheoremstyle{lem}{6pt}{6pt}{\rm}{}{\sffamily}{ }{ }{}
 \theoremstyle{lem}

\newtheoremstyle{case}{6pt}{6pt}{\rm}{}{}{. }{ }{}
 \theoremstyle{case}

 \newtheoremstyle{statement}{6pt}{6pt}{\rm}{}{\sffamily}{ }{ }{}
\theoremstyle{statement}

 \newtheoremstyle{corollary}{6pt}{6pt}{\rm}{}{\sffamily}{ }{ }{}
 \theoremstyle{corollary}

  \newtheoremstyle{defi}{6pt}{6pt}{\rm}{}{\sffamily}{ }{ }{}
 \theoremstyle{defi}
 \newtheorem{defi}{\sc Definition}
 
  \newtheoremstyle{cor}{6pt}{6pt}{\rm}{}{\sffamily}{ }{ }{}
 \theoremstyle{cor}
\newtheorem{cor}{\sc Corollary}[section]

\newtheoremstyle{example}{6pt}{6pt}{\rm}{}{\sffamily}{ }{ }{}
\theoremstyle{example}

\newtheoremstyle{remark}{6pt}{6pt}{\rm}{}{\sffamily}{ }{ }{}
\theoremstyle{remark}

\newtheoremstyle{approximation}{6pt}{6pt}{\rm}{}{\sffamily}{ }{ }{}
\theoremstyle{approximation}

\newtheoremstyle{scheme}{6pt}{6pt}{\rm}{}{\sffamily}{ }{ }{}
\theoremstyle{scheme}

\newtheoremstyle{Algorithm}{6pt}{6pt}{\rm}{}{\sffamily}{ }{ }{}
\theoremstyle{Algorithm}

 \newtheoremstyle{Remark}{6pt}{6pt}{\rm}{}{\sffamily}{ }{ }{}
 \theoremstyle{Remark}

\newtheoremstyle{Lemma}{6pt}{6pt}{\rm}{}{\sffamily}{ }{ }{}
\theoremstyle{Lemma}

\newtheoremstyle{Assumption}{6pt}{6pt}{\rm}{}{\sffamily}{ }{ }{}
\theoremstyle{Assumption}

\newtheoremstyle{Proposition}{6pt}{6pt}{\rm}{}{\sffamily}{ }{ }{}
\theoremstyle{Proposition}

\newtheoremstyle{prop}{6pt}{6pt}{\rm}{}{\sffamily}{ }{ }{}
\theoremstyle{prop}

\newtheoremstyle{rem}{6pt}{6pt}{\rm}{}{\sffamily}{ }{ }{}
 \theoremstyle{rem}

\newtheoremstyle{hypo}{6pt}{6pt}{\rm}{}{\sffamily}{ }{ }{}
 \theoremstyle{hypo}

  \newtheoremstyle{Step}{6pt}{6pt}{\rm}{}{}{ }{ }{}
 \theoremstyle{Step}

 \newtheoremstyle{lema}{6pt}{6pt}{\rm}{}{\sffamily}{ }{ }{}
 \theoremstyle{lema}

\renewcommand{\theequation}{\thesection.\arabic{equation}}
\numberwithin{equation}{section}

\begin{document}


\title{Weak Convergence in the Prokhorov Metric
 of Methods for  \\ Stochastic Differential 
Equations}

\author{{\sc Benoit Charbonneau$^1$, Yuriy Svyrydov$^2$, P. F. Tupper$^3$}\\[2pt]
$^1$Mathematics Department, Duke University, Durham, NC, U.S.A.\\
$^2$Technische Universit\"at M\"unchen, Munich, Germany.\\
$^3$Department of Mathematics, Simon Fraser University, Burnaby, BC, Canada.}

\pagestyle{headings}

\maketitle

\begin{abstract}
{
We consider the weak convergence of numerical methods for stochastic differential equations (SDEs). 
Weak convergence is usually expressed in terms
of the convergence of expected values of  test functions of the
trajectories.  Here we present an alternative formulation of
weak convergence in terms of the well-known Prokhorov metric
 on spaces of random variables.  For a general
class of methods, we establish bounds on the rates of
convergence in terms of the Prokhorov metric. 
In doing so, we revisit the original proofs of weak
convergence and show explicitly how the bounds on the error depend on the
smoothness of the test functions. 
As an application of our result, we  use the Strassen--Dudley theorem
to show that the numerical approximation and the true solution to the system of SDEs
can be re-embedded in a
probability space in such a way that the
method converges there  in a strong sense.    One  corollary of this last result is that the method converges in the Wasserstein distance, another metric on spaces of random variables.  Another corollary establishes rates of convergence for expected values of  test functions assuming only local Lipschitz continuity. 
We conclude with a review of the existing results for pathwise convergence of weakly converging methods
and the corresponding strong results available under re-embedding.}
\end{abstract}

\section{Introduction}

Consider the following system of Ito stochastic differential equations (SDEs)
\begin{equation} \label{eqn:sde0}
dX = a(X) dt + \sum_{r=1}^q \sigma_r(X) dW_r(t), \ \ \ \ X(0)=x_0,
\end{equation}
for $X(t) \in \mathbb{R}^n$,
where the $W_r(t)$ are independent standard Wiener processes.
The simplest numerical method for obtaining approximate solutions to
this system is the Euler--Maruyama method: for $k \geq 0$, timestep $\Delta t$, and 
$\Delta_k W_r = W_r((k+1)\Delta t) - W_r(k \Delta t)$, 
\begin{equation} \label{eqn:strongmethod}
X_{k+1} = X_k + a(X_k) \Delta t + \sum_{r=1}^q \sigma_r (X_k) \Delta_k W_{r},  \ \ \ \ X_0=x_0.
\end{equation}
  For each $k$, $X_k$ is an approximation to $X(k \Delta t)$. The
Euler--Maruyama method converges in the \emph{strong} sense
 because for each
realization of the $W_r(t)$, the method gives an approximation to
the exact solution of the SDE with that same realization.  In
particular, as shown in \cite[p.\ 342]{kloeden},
\begin{equation} \label{eqn:meansq}
(\mathbb{E} (X(T) - X_{T/\Delta t})^2)^{1/2} \leq C \Delta t^{1/2},
\end{equation}
under certain assumptions on the coefficients $a$ and $\sigma$.
In order for such a result to be possible $X(T)$ and $X_{T/\Delta t}$ must be defined on the same probability space.

Another way to quantify convergence of a numerical method is to consider the distribution of the random variable generated by the numerical method and see how close it is to the distribution of the  true trajectory at the corresponding point in time.  This concept is known as \emph{weak} convergence.  The typical way to quantify this is through test functions.  For example, for sufficiently smooth functions $f$ the Euler--Maruyama method (\ref{eqn:strongmethod}) satisfies 
\begin{equation} \label{eqn:weakconv}
| \mathbb{E} f(X_{T/\Delta t}) - \mathbb{E} f(X(T)) | \leq C_f \Delta t
\end{equation}
for some constant $C_f$ depending on $f$.  See  \cite[p.\ 473]{kloeden}, or \cite{talay} for an important earlier reference.

Strong convergence of a method implies weak convergence, but the converse is not  true.  For example, let $N_{rk}$  be independent identically distributed random variables with $N_{rk} = \pm 1$ with probability $1/2$, and set
\begin{equation} \label{eqn:weakmethod} 
X_{k+1} = X_k + a(X_k) \Delta t + \sum_{r=1}^q \sigma_r (X_k) \Delta
t^{1/2} N_{rk},
\end{equation}
We can define the $N_{rk}$ in various ways.   One possibility is for them to be independent of  the Wiener processes driving the SDE (\ref{eqn:sde0});  another is to choose them to be 
$N_{rk} = \mbox{sgn}( \Delta_k W_k)$.  In either case,
this method, which we call weak Euler--Maruyama, does not converge strongly to the solution of the SDE.  
 However, (\ref{eqn:weakconv}) still holds whatever the relation between the $N_{rk}$ and the Wiener processes, even if they are defined on different probability spaces.

In this paper, we present a different formulation of weak convergence of numerical method in terms of the Prokhorov metric \cite{billingsley}.  For any two random elements of $\mathbb{R}^n$ the Prokhorov metric gives a quantitative measurement of how far apart their distributions are.  Its importance in probability theory \cite{dudley} is that convergence in the Prokhorov metric is equivalent to convergence in distribution (or weak convergence, as it is sometimes known).  
Before we discuss our results
 we first review some of the definitions and  facts of convergence in distribution in metric spaces.  See Billingsley's book \cite{billingsley} for details.  

Consider a metric space $S$ with metric $d$, such as $\mathbb{R}^n$ with the Euclidean metric.   We say that a sequence of random elements $X_n$ in $S$ \emph{converges in distribution} 
to $X$ in $S$ if for all bounded continuous $f : S \rightarrow \mathbb{R}$
\begin{equation} \label{eqn:def1}
\mathbb{E}f(X_n) \rightarrow \mathbb{E} f(X),
\end{equation}
as $n \rightarrow \infty$.  An equivalent definition of convergence in distribution of $X_n$ to $X$ is that for all Borel sets $A$ of $S$ with $\mathbb{P}(X \in \partial A) = 0$ we have
\begin{equation} \label{eqn:def2}
\mathbb{P}(X_n \in A) \rightarrow  \mathbb{P}(X \in A)
\end{equation}
as $n \rightarrow \infty$.  (Here $\partial A$ is the boundary of the set $A$.)  The assumption of boundedness on $f$ may seem excessive, but in the presence of uniform bounds on the moments of $X_n$ and $X$, convergence in distribution implies that (\ref{eqn:def1}) holds for more general continuous $f$  \cite[p.\ 86]{babydurrett}.

 It is not obvious from either of the above  definitions of weak convergence how to measure the speed with which a sequence $X_n$ converges in distribution to $X$, since the rate at which limits (\ref{eqn:def1}) and (\ref{eqn:def2}) occur depends on $f$ and $A$ respectively.
The Prokhorov metric is one  way to define the distance between the distributions of two random elements, and thus allows us to quantify convergence in distribution.
For any two random elements $X$ and $Y$ of $S$ let $\rho(X,Y)$ be the Prokhorov distance between them  (see Section 2 for the definition).  This distance is zero if and only if $X$ and $Y$ have the same distribution, that is if $\mathbb{P}(X \in A) = \mathbb{P}(Y \in A)$ for all Borel sets $A$.  Moreover, if $X_n$ is a sequence of random elements in a separable metric space
$S$, $\rho(X_n,X) \rightarrow 0$ if and only if
$X_n$ converges in distribution to $X$.
Thus we say that the Prokhorov metric metrizes convergence in distribution.

In our case,
we view the solution of the system of SDEs at time $T$ as random vector (a random element of 
$\mathbb{R}^n$),  and likewise for the numerical solution at time $T$.  Then we ask how the Prokhorov distance between  $X(T)$ and $X_{T/\Delta t}$ depends on $\Delta t$.  Our main result in Section 3 shows that  
 the usual definition of weak convergence in terms of test functions implies convergence in the Prokhorov metric, and we provide a bound on the rate.  One important component of our proof is determining how exactly the constant $C_f$ in (\ref{eqn:weakconv}) depends on $f$ in the usual proofs of weak convergence \cite{milstein,kloeden}.

In Section 4 we show one consequence of our main result concerning re-embedding trajectories of the SDEs and of the numerical method in a new probability space.  Two random vectors (such as $X(T)$ and $X_{T/\Delta t}$) may either not be close to each other on a realization-by-realization basis, or may be defined on completely different probability spaces.
However, it is possible to define new random vectors $Y$ and $Z$ jointly on a new probability space such that $Y$ has the same distribution as $X(T)$ and $Z$ has the same distribution as $X_{T/\Delta t}$.  This construction is called a \emph{re-embedding} of $X(T)$ and $X_{T/\Delta t}$ in a new probability space.
After re-embedding the random vectors may be close together in a strong sense and we can look at how
quantities like  $\mathbb{E} |Y-Z|$ or $\mathbb{P}(|Y-Z|> \alpha)$ behave as $\Delta t$ varies. 
The Strassen--Dudley theorem says that if two random variables are close in the Prokhorov metric, then there is a re-embedding of them into another probability space for which they are close in probability.  A bound on some higher moment of $Y$ and $Z$ then gives that $\mathbb{E}|Y-Z|$ is small.  
Using our bound in the Prokhorov metric and the Strassen-Dudley theorem, we show that a method with the usual weak convergence of order $p$ converges strongly after re-embedding with order $\frac{p}{2p+3}- \epsilon$ for any $\epsilon>0$.  This is equivalent to proving a rate of convergence in the Wasserstein distance (see Section 4 for a definition).
We also use re-embedding to  establish rates for the  convergence of expectations of test functions  requiring only  local Lipschitz continuity and polynomial growth.

Finally, in Section 5, we discuss the corresponding result for weak convergence of entire numerical  trajectories on $[0,T]$  to exact trajectories of the original system.  Convergence in distribution follows directly from a result of Stroock and Varadhan, which we review.  However, to the best of our knowledge there is no bound available for this rate of convergence for general weakly convergent methods.  (Though see \cite{gladyshev,kanagawa} for some results for the strong Euler--Maruyama method.)
Applying Skorohod's theorem  gives the corresponding strong convergence result for the trajectories embedded in another probability space, though again without a rate.

\section{Metrics on Spaces of Random Elements}

Consider a metric space $(S,d)$ with metric $d$.  A random element of $S$ is a measurable function $X: \Omega \rightarrow S$ where $(\Omega,\mathcal{F},\mathbb{P})$ is some probability space.  For example, if $S$ is $\mathbb{R}^n$ with the metric $d(x,y) =|x-y|$, then random elements $X$ are called random vectors.  Even if two random elements $X$ and $Y$ of a metric space $S$
are not close on a realization-by-realization basis, we may still wish to compare their distributions.  So we define a metric  on the space of random elements of $S$.
Note that  there are two distinct metrics involved: $d$ which is a metric on the original space $S$ and 
 another  which is  a metric on  the space of  random elements of $S$. 
In this section, we  first define the well-known Prokhorov metric $\rho$, which is defined for any underlying metric space. 
Then we introduce the metrics $\beta_l$ for non-negative integers $l$, when the underlying metric space is $\mathbb{R}^n$.
The latter are similar to the metric $\beta$ of Fortet and Mourier \cite{fortet}; see \cite[Sec. 11.3]{dudley}.

For a set $A \subset S$ we define $A^{\epsilon}$, $\epsilon>0$, the set of all points within distance $\epsilon$ of $A$ by 
\begin{equation}\label{eqn:epdef}
A^{\epsilon} = \{ x \in S \mid \inf_{y \in A}d(x,y) \leq \epsilon \}.
\end{equation}
The Prokhorov metric is defined as follows.

\begin{defi} \label{def:prokhorov} For random variables $X$ and $Y$ in $S$
\[
\rho(X,Y) : = \inf \{ \epsilon \mid \mathbb{P} (X \in A) \leq \mathbb{P}(Y \in A^{\epsilon}) + \epsilon, \mbox{ for all $A$ closed} \}.
\]
\end{defi}

If we identify random elements of $S$ that have the same distribution, then $\rho$ is a metric on the set of random elements \cite[p.\ 394]{dudley}.    If $(S,d)$ is separable (as are all examples in this paper) random elements $X_n$ converge in distribution to $X$ if and only if $\rho(X_n,X) \rightarrow 0$ \cite[p.\ 395]{dudley}.  
 Note that $\rho(X,Y) \leq 1$ always.

Here is the Strassen--Dudley Theorem as proven in \cite{dudley,strassen},  used later in this section and in Section~\ref{sec:strongcon}.

\begin{thm} \label{thm:strass} (\cite[p.\ 73]{billingsley})  Let $(S,d)$ be a separable metric space.
 If $X$ and $\bar{X}$ are random elements of $S$ with 
  $\rho(X,\bar{X}) < \alpha$,
then there are random elements $Y$ and $Z$ of $S$ defined on a common probability space 
  such that  $Y$ has the same distribution as $X$, $Z$ has the same distribution as $\bar{X}$ and 
\[
\mathbb{P} ( |Y-Z| > \alpha ) < \alpha. 
\]
\hfill $\square$
\end{thm}

We now define a class of metrics $\beta_l$ on random vectors, that is, random elements of the metric space  $(\mathbb{R}^n,| \cdot |)$.
 Let $f \colon \mathbb{R}^{k} \rightarrow \mathbb{R}$.  Let $\alpha$ be a vector of length $k$ with non-negative integer components.  Let $| \alpha |:=\sum_i \alpha_i$ and 
\[
D^{\alpha} f := \frac{ \partial^{|\alpha|} f}{\partial x_1^{\alpha_1} \ldots
\partial x_n^{\alpha_n}}.
\]
If we wish to emphasize the argument of $f$ in our notation we use $D_x^\alpha$ instead of $D^\alpha$.
For $l \geq 0$
and $f \colon \mathbb{R}^n \rightarrow  \mathbb{R}$ let
\begin{equation} \label{eqn:normdef}
\| f \|_l := \sum_{|\alpha| \leq l} \sup_{x \in \mathbb{R}^n} |  D^\alpha  f (x) |.
\end{equation}
\begin{defi} For random vectors $X$ and $Y$ in $\mathbb{R}^n$
  and for $l \geq 0$ we let
\[
\beta_l(X,Y) = \sup_{\| f\|_l \leq 1}  \left|
\mathbb{E} f(X) - \mathbb{E} f(Y) \right|
\]
\end{defi}
It is straightforward to check that $\beta_l$ is a metric on the space of random variables.

The following theorem is the main result in this section, and allows us to show in the next section that solutions generated by weak numerical methods converge in the Prokhorov metric.  

\begin{thm} \label{thm:produd} For each $l\geq 0$ there is a constant $C>0$ such that
for any random vectors $X$, $Y$ in $\mathbb{R}^n$,
\[
\rho(X,Y) \leq C \beta_l (X,Y)^{1/(1+l)}.
\]
\end{thm}

\begin{proof}  Here we closely follow \cite[p.\ 396]{dudley}.  
Consider any closed set $K$ in $\mathbb{R}^n$ and $\epsilon \in (0,1]$.   Let $K^\epsilon$ be defined as in Equation (\ref{eqn:epdef}).
 The lemma following this theorem shows that there is a smooth function $f$  and a constant $C$ such that  depends on $n$ but not on $\epsilon$ or $K$ such that 
\[
\mathbf{1}_K(x) \leq f(x)  \leq \mathbf{1}_{K^\epsilon}(x)
 \ \ \ \ \   \mbox{ and }  \ \ \ \ \ \ 
\|f\|_l \leq C \epsilon^{-l}.
\]
Without loss of generality we assume that $C\geq 1$.
We now use the function $f$ to establish the required bound.
  For any random variables $X$ and $Y$
\begin{align*}
\mathbb{P}( Y \in K ) & \leq  \mathbb{E} f(Y) \\
& \leq  \mathbb{E} f(X) + \|f\|_l \beta_l (X,Y) \\
& \leq  \mathbb{P}( X \in K^\epsilon) + C \epsilon^{-l} \beta_l(X,Y).
\end{align*}
 So for any $\epsilon \in (0,1]$
 \begin{equation}\label{eqn:helper}
 \rho(X,Y) \leq  \max(\epsilon, C \epsilon^{-l} \beta_l(X,Y)). 
\end{equation}
  Now if $\beta_l(X,Y) >1$, since $C \geq 1$ and $\rho(X,Y)\leq 1$ the result is immediately true.
So we assume that  $\beta_l(X,Y) \leq 1$ and choose $\epsilon$ so that $\beta_l(X,Y) = \epsilon^{l+1}$.  Then $\epsilon \leq 1$ and (\ref{eqn:helper}) gives us 
 $ \rho(X,Y) \leq \epsilon \max(1,C )$.  So $\rho(X,Y) \leq C \beta_l (X,Y)^{1/(1+l)}$ as required. 
\end{proof}

\begin{lemma}  For each closed set $K \subset \mathbb{R}^n$ there is a 
parametrized family of functions $f_\epsilon(x)$ for $\epsilon \in (0,1]$ such that 
\begin{equation} \label{eqn:sandwich}
\mathbf{1}_K(x) \leq f_\epsilon (x)  \leq \mathbf{1}_{K^\epsilon}(x),
\end{equation}
and there is a constant $C$ depending on $n$ but not on $\epsilon$, $K$, or $l$ such that
\begin{equation} \label{eqn:erf}
\|f_\epsilon \|_l \leq C \epsilon^{-l}.
\end{equation}
\end{lemma}
\begin{proof}
We use the method of mollifiers; see, for example \cite[p.\ 629]{evans}.
Define $\eta \colon \mathbb{R}^n \rightarrow \mathbb{R}$ by
\[
\eta(x) := \left\{ \begin{array} {lcc}
D \exp \left( \frac{1}{\|x\|^2-1} \right) & \mbox{if} & \|x\|<1 \\
0 & \mbox{if} & \|x\|\geq 1,
\end{array} \right.
\]
where $D>0$ is selected so that $\int_{\mathbb{R}^n} \eta(x) dx = 1$. 
The mollifier $\eta \in C^\infty$ is positive with support in the unit ball about the origin.
Define
\[
\eta_{\epsilon}(x) := \frac{1}{\epsilon^n} \eta\left(\frac{x}{\epsilon}\right).
\]
This function is in $C^{\infty}$, has support on the ball of radius $\epsilon$ about the origin, and also has integral $1$.

Let $K'$ be the closure of $K^{\epsilon/2}$ and  let 
\[
f_\epsilon(x) : = \int_{\mathbb{R}^n} \eta_{\epsilon/2} (y)
\mathbf{1}_{K'}(x-y) dy.
\]
The function $f_\epsilon$ is 1 on $K$, 0 on $\mathbb{R}^n \setminus
K^{\epsilon}$, and between zero and one elsewhere.  So $f_\epsilon$ satisfies the condition of Equation
(\ref{eqn:sandwich}). 

In \cite[p.\ 630]{evans} it is shown that 
\[
D^\alpha f_\epsilon (x) = \int_{\mathbb{R}^n} D_x^\alpha \eta_{\epsilon/2}(x-y) \mathbf{1}_{K'}(y) dy.
\]
So
\begin{align*}
|D^\alpha f_\epsilon (x) |  & \leq 
 \int_{\mathbb{R}^n} | D^\alpha_x \eta_{\epsilon/2}(x-y) | dy   
    =  \int_{\mathbb{R}^n} (\epsilon/2)^{-n} 
    | D^\alpha_x \eta( 2(x-y) /\epsilon )| dy \\
 &  =  (\epsilon/2)^{-|\alpha|}  \int_{\mathbb{R}^n} | D^{\alpha}_z \eta(z)| dz, 
 \end{align*}
 where we have used the change of variables $z=2(y-x)/\epsilon$,  $D^{\alpha}_z = (-\epsilon/2)^{|\alpha|} D^{\alpha}_x$.  The integral in the last expression is finite and does not depend on $K$ or $\epsilon$.
 Summing over all $\alpha$ with $|\alpha| \leq l$ gives us
 \[
 \| f \|_l \leq \sum_{|\alpha |\leq l} C_{\alpha} \epsilon^{-|\alpha|} \leq C \epsilon^{-l}
 \] 
 for some constants $C_\alpha, C$, for all $\epsilon \in (0,1]$.
 \end{proof}

For completeness we include the following theorem which together with Theorem~\ref{thm:produd} shows that the metrics $\rho$ and $\beta_l$ induce the same topology on the space of random elements of $\mathbb{R}^n$.  Thus, as for $\rho$, $\beta_l(X_n,X) \rightarrow 0$ if and only if $X_n$ converges to $X$ in distribution.  This result is analogous to Theorem 11.6.5 in \cite{dudley}.

\begin{thm}
For all $l \geq 1$, and random $X$ and $\bar{X}$ in $\mathbb{R}^n$, the metrics $\rho$ and $\beta_l$ satisfy
\[
\beta_{l} (X,\bar{X}) \leq  2 \rho(X,\bar{X}).
\]
\end{thm} 
\begin{proof}   
Let $\rho(X,\bar{X})= \epsilon$.  Using Theorem \ref{thm:strass}, let $Y$ and $Z$ be random vectors on the same probability space with the same distributions as $X$ and $\bar{X}$ respectively
such that
\(
\mathbb{P}( |Y-Z|>\epsilon)<\epsilon.
\)
Then 
\begin{align*}
\beta_l(X,\bar{X}) & \leq  \sup_{\|f\|_l \leq 1} \mathbb{E} | f(Y) - f(Z)|  \\
& =  \sup_{\|f\|_l \leq 1}  \left\{ \mathbb{E} \left[  \mathbf{1}_{|Y-Z|>\epsilon} |f(Y)-f(Z)| \right]+ 
        \mathbb{E} \left[ \mathbf{1}_{|Y-Z| \leq \epsilon} |f(Y)-f(Z)| \right] \right\}  \\
& \leq  \sup_{\|f\|_l \leq 1} \left\{
\mathbb{E} \left[  \mathbf{1}_{|Y-Z|>\epsilon}  \right]  2 \sup_x |f(x)| + \epsilon   \max_{y,z} \frac{|f(y)-f(z)|}{\|y-z\|} \right\} \\
& \leq   \sup_{\|f\|_l \leq 1} \left\{
2 \epsilon  \sup_x |f(x)| + \epsilon  \sum_{i} \sup_x |D^{i}(f(x))| \right\} \\
& \leq  2 \epsilon, 
\end{align*}
as required. \end{proof}
 
\section{Convergence of Numerical Methods}

Here we prove our result on the convergence in the Prokhorov metric of numerical approximations to exact solutions of SDEs. 
 We consider the system of Ito SDEs
\begin{equation}\label{eqn:sde}
dX = a(X) dt + \sum_{r=1}^q \sigma_r(X) dW_r(t),
\end{equation}
where $X(t) \in \mathbb{R}^n$, $a \colon \mathbb{R}^n \rightarrow \mathbb{R}^n$, $\sigma_r \colon \mathbb{R}^n \rightarrow \mathbb{R}^{n \times n}$ for all $r$. The  $W_r$, $r=1,\ldots,q$  are mutually independent  standard Wiener processes.  We set the initial condition to be $X(0)=x_0 \in \mathbb{R}^n$.

To prove our convergence theorem 
we build on a weak convergence result from \cite{milstein}.
This result is expressed for a rather general method for the system (\ref{eqn:sde}):
\begin{equation} \label{eqn:method}
X_{k+1} = X_k + \bar{\delta}(X_k, \Delta t ; \xi_k),
\end{equation}
with $X_0=x_0$.
Here $\bar{\delta}$ is vector-valued function and $\xi_k$, $k=0,1, \ldots$  is a sequence of independent random vectors.  Usually we suppress the $\xi_k$ from the notation and view $\bar{\delta}(X_k,\Delta t)$ as a random vector.  We denote its $i$th component by $\bar{\delta}_i(X_k,\Delta t)$.
Here $X_k$ is intended to be an approximation to $X(k \Delta t)$.
In the following we use $\delta$ to denote the increment of the true solution over a time interval:
for the solution $X$ to Equation (\ref{eqn:sde}) with $X(0)=x$, set
\[ 
\delta (x,\Delta t) = X(\Delta t)-X(0).
\] 
Thus $\delta(x,\Delta t)$, like $\bar{\delta}(x,\Delta t)$, is a random vector.  The $i$th component of $\delta(x,\Delta t)$ is denoted $\delta_i(x,\Delta t)$.

Theorem~\ref{thm:milstein} below gives a rate of convergence of $\mathbb{E}f(X_k)$ to $\mathbb{E}f(X(k \Delta t))$ in which the dependence of the constant on $f$ is given.   This result is an corollary of the result  of \cite[p.\ 100]{milstein} or \cite[p.\ 473]{kloeden} in which the dependence of the constant on $f$ is not made explicit.  Here, by making stronger assumptions on the coefficients $a$ and $\sigma_r$, we show that the constant is linear in  $\|f\|_{2p+2}$ where $p$ is the order of the method.  (See Section 2 for a definition of $\| \cdot \|_{2p+2}$.)

\begin{thm} \label{thm:milstein}
Let $T>0$ be fixed.
Suppose that  \\
(a) the coefficients $a$ and $\sigma_r$ of  the system of SDEs (\ref{eqn:sde})
have globally Lipschitz derivatives up to and including order $2p+2$;\\
(b) there is some scalar function $K(x)$ with at most polynomial growth as $x \rightarrow \infty$
such that 
\[
\left| \mathbb{E} \left( \prod_{j=1}^s \delta_{i_j}(x,\Delta t) - \prod_{j=1}^s \bar{\delta}_{i_j}(x,\Delta t) \right) \right| \leq K(x) \Delta t^{p+1},
\]
for $s=1, \ldots, 2p+1$ and 
\[
\mathbb{E} \prod_{j=1}^{2p+2} | \bar{\delta}_{i_j}(x,\Delta t) | \leq K(x) \Delta t^{p+1};
\]
(c) for all $m \geq 1$ the expectations $\mathbb{E} | X_k |^{2m}$ exist and are uniformly bounded with respect to $\Delta t$ and $k =0,1,\ldots,\lfloor T/\Delta t \rfloor.$ \\
(d) the function $f(x)$ together with its partial derivatives of order up to and including $2p+2$ are bounded. 
Then for  all $k\Delta t \in [0,T]$
\[
| \mathbb{E} f(X(k \Delta t)) - \mathbb{E} f( X_k ) | \leq C \|f\|_{2p+2}  \Delta t^p.
\]
The constant  $C$ depends on $x$, $a$, $\sigma_r$ and $T$ but not on $f$ and $\Delta t$.  \end{thm}

\begin{proof}  
We define  $Y(x,t)$ to be $X(k \Delta t)$ where $X$ is the solution of (\ref{eqn:sde}) with initial condition $X(t)=x$, $t\leq k \Delta t $.
Define the function $u(x,t)$ by 
\[
u(x,t) := \mathbb{E} f (Y(x,t)).
\]
If follows from the proof of Theorem 2.1 in \cite[p.\ 100]{milstein} that
\[
| \mathbb{E} f(X(k \Delta t)) - \mathbb{E} f(X_k) |
 \leq A \Delta t^p \max_{t \in [0,k \Delta t]}\| u(\cdot,t) \|_{2p+2}
\]
for some constant $A$ not depending on $f$.

In \cite[p.\ 223]{kunita} it is shown that if the coefficients of $a$ and $\sigma_r$ of the system of SDEs (\ref{eqn:sde}) have globally Lipschitz continuous derivatives up to order $2p+2$ (condition (a))\ then 
 $Y(t,x)$ has continuous derivatives with respect to $x$ up to order $2p+2$, almost surely.  
 
 Let $\partial_i$ denote differentiation of a function with respect to its $i$th argument.
 Formally, we can differentiate $u$ with respect to $x$ to obtain
 \begin{align*}
 \partial_i u & =    \sum_a \mathbb{E} \left[ \bigl( \partial_a f(Y) \bigl)(\partial_i Y_a) \right], \\
\partial_i \partial_j  u & =   \sum_{a,b} \mathbb{E}  \left[ \bigl(\partial_a \partial_b f (Y) \bigl)(\partial_i Y_a)(\partial_j Y_b) \right]+   \sum_{a} \mathbb{E} \left[ \bigl( \partial_a f(Y)  \bigl)(\partial_i \partial_j Y_a) \right],
 \end{align*} 
 and so forth, using the product and chain rules.
To justify the formal differentiations, we need only observe that all multi-derivatives of $f$ up to order $2p+2$ are bounded, and remark that it follows from \cite{kunita} that all moments of the derivatives of $Y$ up to order $2p+2$ are finite.  The exchange of differentiation with expectation then follows in each case by Fubini's theorem \cite[p.\ 222]{williams}.
Applying the Cauchy--Schwarz inequality to each term gives 
\[
\sup_{x \in \mathbb{R}^n} | D^{\alpha} u | \leq \sum_i \left\{ \mathbb{E} (D^{\beta_i} f)^2  \right\}^{1/2} E_{\beta_i},
\]
where $\beta_i$ are a sequence of multi-indices with $| \beta_i | \leq |\alpha|$ and 
 $E_{\beta_i}$ are some constants independent of $f$.
So 
\[
\sup_{x \in \mathbb{R}} | D^{\alpha} u |  \leq F_\alpha \| f \|_{|\alpha|},
\]
for some constants $F_\alpha$ independent of $f$.
Summing this inequality over all $\alpha$ with $|\alpha| \leq 2p+2$ gives us the result.
\end{proof}

Putting Theorems \ref{thm:produd} and \ref{thm:milstein} together gives us our conclusion for this section.

\begin{cor}  \label{thm:fdependence}
Let conditions (a), (b), (c) of Theorem \ref{thm:milstein}  be satisfied.  Then for some constant $K$
\[
\rho(X(k \Delta t),X_k) \leq K \Delta t^{p/(2p+3)},
\]
for all $k \Delta t \leq T$.
\end{cor}
\begin{proof} By Theorem~\ref{thm:milstein} and the definition of $\beta_{l}$, we have that
\[
\beta_{2p+2}(X(k \Delta t), X_k) \leq C \Delta t^p.
\]
Applying Theorem \ref{thm:produd} with $l=2p+2$ then gives the result. \end{proof}

\section{Strong Convergence} \label{sec:strongcon}

In this section, we apply the Strassen--Dudley theorem to show that,  after being re-embedded in another probability space, weakly converging methods for stochastic differential equations
converge strongly with a reduced order. 
This re-embedding 
immediately gives a rate of convergence in the Wasserstein distance.
As a corollary, we establish a rate of convergence of $\mathbb{E}
f(X_{T/\Delta t})$ to $\mathbb{E} f(X(T))$ that requires only that $f$ is locally Lipschitz with a polynomial growth condition.

\begin{thm} \label{thm:easy} Let conditions (a), (b), (c) of Theorem \ref{thm:milstein} be satisfied.   There is a probability space on which random vectors $Y$ and $Z$ are defined such that $Y$ has the same distribution as $X(T)$ and $Z$ has the same distribution as $X_{T/\Delta t}$ and, for any $\epsilon >0$
\[
\mathbb{E} |Y - Z| \leq C \Delta t^{ \frac{p}{(2p+3)} - \epsilon},
\]
for some constant $C$, for all sufficiently small $\Delta t$.
\end{thm}
\begin{proof}  Let $\alpha=K \Delta t^{p/(2p+3)}$ where $K$ is as in Corollary~\ref{thm:fdependence}.  Theorem \ref{thm:strass} together with Corollary~\ref{thm:fdependence} establish the existence of the random vectors $Y$ and $Z$ with the correct distributions such that
\[
\mathbb{E}( \mathbf{1}_{|Y-Z|>\alpha}) = \mathbb{P}(|Y-Z| > \alpha)< \alpha.
\]
Choose $\epsilon>0$.
Now, the conditions of Theorem \ref{thm:milstein} ensure that both $Y$ and $Z$ and hence $|Y-Z|$ have finite  moments of all orders independent of $\Delta t$.  Choose real numbers $q_1,q_2>1$ such that $1/q_1 + 1/q_2 = 1$ and $\frac{p}{(2p+3)} \frac{1}{q_2} \geq \frac{p}{(2p+3)} -\epsilon$. 
Using H\"older's inequality, we obtain
\begin{align*}
\mathbb{E} |Y - Z|  & =   \mathbb{E} \left[ |Y-Z| \mathbf{1}_{|Y-Z| > \alpha} \right] +
 \mathbb{E} \left[ |Y-Z| \mathbf{1}_{|Y-Z| \leq \alpha}  \right] \\
& \leq   (\mathbb{E} |Y-Z|^{q_1})^{1/q_1} ( \mathbb{E} \mathbf{1}_{|Y-Z|>\alpha})^{1/q_2} + \alpha \\
& \leq  (\mathbb{E} |Y-Z|^{q_1})^{1/q_1}  \alpha^{1/q_2} + \alpha \\
& \leq  E \alpha^{1/q_2} = D \Delta t^{\frac{p}{(2p+3)} \frac{1}{q_2}}
\leq  C \Delta t^{\frac{p}{(2p+3)}-\epsilon},
\end{align*}
for all sufficiently small $\Delta t$,
as required.
\end{proof}

Applying this theorem to the case of weak Euler--Maruyama  (see Equation (\ref{eqn:weakmethod}))
with $p=1$ implies a strong rate of convergence after re-embedding of 
 $1/5 - \epsilon$  for any $\epsilon>0$.  

Now we express our result in terms of the Wasserstein distance, also known as the Wasserstein-1 distance, the Monge--Wasserstein distance \cite[p.\ 420]{dudley}, or the Kantorovich--Rubinstein
distance \cite[p.\ 206]{villani}.
To define this metric, let $X$ and $Y$ be random elements of $S$ and let $M(X,Y)$ be the set of all probability measures $\mu$ on $S \times S$ such that 
the marginals of $\mu$ are the probability measures induced on $S$ by $X$ and $Y$ respectively.
 Then
the Wasserstein distance is
\[
W(X,Y)= \inf_{\mu \in M(X,Y)}  \mathbb{E}_{\mu} d(x,y)
= \inf_{\mu \in M(X,Y)}  \int d(x,y)\, d\mu(x,y).
\]
where $\mathbb{E}_{\mu}$ denotes expectation with respect to the measure $\mu$ for  $(x,y) \in S \times S$.
In words, the Wasserstein distance is the minimal $L^1$ distance between $X$ and $Y$ after re-embedding.  Therefore Theorem \ref{thm:easy} shows the following:
\begin{cor} Let conditions (a), (b), (c) of Theorem \ref{thm:milstein} be satisfied.  Then  for any $\epsilon>0$
\[
W(X(T),X_{T/\Delta t}) \leq C \Delta t^{\frac{p}{2p+3}-\epsilon},
\]
for some constant $C$, for sufficiently small $\Delta t$.  \hfill $\square$ \end{cor} 

As another corollary to Theorem \ref{thm:easy}, we show that  $\mathbb{E} f (X_{T/\Delta t}) \rightarrow \mathbb{E} f (X(T))$ given some polynomial growth conditions on $f$, even when $f$ is only locally Lipschitz.  This result  is like the usual weak convergence result \cite[p.\ 100]{milstein}, but with a relaxed  smoothness requirement on $f$ and a reduced rate.  Compare with Mikulevicius and Platen's result \cite[p.\ 460]{kloeden}  or Bally and Talay's result  \cite{ballytalay}, both of which only apply to strong Euler--Maruyama (see Equation (\ref{eqn:strongmethod})).

\begin{cor} Let conditions (a), (b), (c) of Theorem \ref{thm:milstein} be satisfied.   Let $f$ be locally Lipschitz with 
\[
| f(x) - f(y) | \leq L_R |x-y|, 
\]
whenever $|x|\leq R$ and $|y| \leq R$,
where
\[
L_R \leq C (1+R^\kappa),
\]
for some constants $C$ and $\kappa$.  
  Then for any $\epsilon>0$
\[
| \mathbb{E} f(X_{T/\Delta t}) - \mathbb{E} f(X(T)) | \leq K \Delta t^{\frac{p}{(2p+3)}-\epsilon},
\]
for some constant $K$, for sufficiently small $\Delta t$.
\end{cor}
\begin{proof}
 Let $\alpha=K \Delta t^{p/(2p+3)}$ where $K$ is as in Corollary~\ref{thm:fdependence}.  
Let $Y$ and $Z$ be as in the proof of Theorem~\ref{thm:easy}, so that $Y$ has the same distribution as $X(T)$, $Z$ has the same distribution as $X_{T/\Delta t}$, and $\mathbb{P}(|Y-Z| > \alpha)< \alpha$.  Let $M = \max(|Y|,|Z|)$. 
We immediately have
\[
| \mathbb{E} f(X_{T/\Delta t}) - \mathbb{E} f(X(T)) | =  | \mathbb{E} f(Y) - \mathbb{E} f(Z) | 
 \leq   \mathbb{E} |f(Y) -  f(Z) |.
 \]
  Let $R$ be an arbitrary radius which we shall fix later.  We can split the quantity of interest into three terms: 
\begin{align*}
 \mathbb{E} |f(Y) -  f(Z) | & \leq 
 \mathbb{E} |f(Y) -  f(Z) | \mathbf{1}_{|Y-Z|\leq \alpha,M \leq R} +
\mathbb{E} |f(Y) -  f(Z) | \mathbf{1}_{|Y-Z| >\alpha,M \leq R} \\
 & \phantom{=}  +\mathbb{E} |f(Y) -  f(Z) | \mathbf{1}_{M>R}  \\
& =: T_1 + T_2 +T_3.
 \end{align*}
To bound the first term, note that
\begin{align*}
T_1 & \leq  L_R \alpha \leq  \alpha C( 1+ R^\kappa).
\end{align*}
To bound the second and third terms, note that  
\begin{align*}
|f(x)|  & \leq |f(0)| + |f(x)-f(0)|  \\
& \leq |f(0)| + L_{|x|} |x| \\
& \leq  D ( 1 + |x|^{\kappa+1}),
\end{align*}
for all $x$ for some constant $D$.  Then
\begin{align*}
T_2 & \leq \mathbb{E} \left\{ |f(Y)| + |f(Z)| \right\} \mathbf{1}_{|Y-Z| >\alpha,M \leq R} \\
& \leq  \mathbb{E} \left\{ D (1+  |Y|^{\kappa+1} ) + D( 1+ |Z|^{\kappa+1}) \right\} \mathbf{1}_{|Y-Z| >\alpha,M \leq R} \\
& \leq \mathbb{E}  2 D (1+  R^{\kappa+1} )  \mathbf{1}_{|Y-Z| >\alpha}  \leq 2 D( 1+ R^{\kappa+1}) \alpha.
\end{align*}
To bound the third term, we use the fact that all moments $Y$ and $Z$ are finite and the bounds on the moments of $Z$ are  independent of $\Delta t$.  For any exponent $q>\kappa + 1$ (which we shall choose later), we let $m_q$ denote this bound so that $\mathbb{E}|Y|^q \leq m_q$ and $\mathbb{E}|Z|^q \leq m_q$.  These inequalities in turn imply that $\mathbb{E}M^q \leq 2 m_q$.
Then
\begin{align*}
T_3 & \leq  \mathbb{E} \left\{| f(Y)| + | f(Z)| \right\} \mathbf{1}_{M>R} \\
&  \leq \mathbb{E} 2 D ( 1+ M^{\kappa+1}) \mathbf{1}_{M > R}\\
&  = 2D \mathbb{E} \frac{( 1+ M^{\kappa+1})}{M^q} M^q \mathbf{1}_{M>R} \\ 
& \leq 2 D \frac{(1+ R^{\kappa+1})}{R^q} 2 m_q \leq C_q  R^{\kappa+1 - q},
\end{align*}
for all sufficiently large $R$.
Putting these three bounds together gives
\begin{align*}
\mathbb{E} | f(Y) - f(Z) |  \leq \alpha E  R^{\kappa+1} + C_q R^{\kappa+1 - q},
\end{align*}
for sufficiently large $R$.
We get to choose both $R$ and $q$.    For any $q \geq \kappa$, if we choose $R = \alpha^{-1/q}$ 
we obtain
\[
\mathbb{E} | f(Y) - f(Z) |  \leq E \alpha^{1 - (\kappa+1)/q} + C_q \alpha^{1- (\kappa+1) / q}.
\]
Choosing $q$ large enough gives the desired result.
\end{proof}

%


\section{Pathwise Convergence in Distribution}

The results in previous sections concern the pointwise weak convergence of numerical methods for SDEs, that is, convergence at each point in time $t$. 
A stronger result is that entire trajectories generated by the numerical method weakly converge to those of the system of SDEs.
Stroock and Varadhan prove pathwise convergence in distribution of numerical methods in great generality in \cite{SandV} but they do not provide a rate.   Here we review their result and apply a re-embedding theorem to establish the corresponding strong result for embedded random paths.  Since no rate appears to be established for Stroock and Varadhan's result, we do not phrase results in terms of the Prokhorov metric and instead just consider convergence in distribution.  Moreover, we can use Skorohod's theorem for re-embedding rather than the Strassen--Dudley theorem.  The latter gives precise rates of strong convergence but the former allows one to construct a whole sequence of random paths and their limit on one probability space.

First we review the definition of convergence in distribution in $C^n[0,T]$, the space of continuous, $\mathbb{R}^n$-valued functions on $[0,T]$, \cite{billingsley}.
 For any fixed $T$ and initial
condition $x_0$ the solution to the system of SDEs (\ref{eqn:sde}) gives a random element of
$C^n[0,T]$ 
which we denote by $X$.  For the same $T$ and
initial conditions the numerical method with step-size $\Delta t$ gives a sequence 
$X_k$, $k=0,1,\ldots$.  We define the linear interpolant $\bar{X}_{\Delta t}$ of the values $X_k$ by
\[
\bar{X}_{\Delta t}(t) = X_{\lfloor t/\Delta t \rfloor } + ( t / \Delta t  - \lfloor t/ \Delta t \rfloor ) 
( X_{\lfloor t/\Delta t \rfloor +1} - X_{\lfloor t/\Delta t \rfloor} ),
\]
for $t\in [0,T]$.  Thus $\bar{X}_{\Delta t}$ is a random element of $C^n[0,T]$.
 If we equip $C^n[0,T]$ with the norm $\| \cdot \|_\infty$
we obtain a metric space with metric
\[
d(x,y) := \|x-y\|_\infty = \sup_{t \in [0,T]} | x(t)-y(t)|,
\]
for $x,y \in C^n[0,T]$. 
We say that $\bar{X}_{\Delta t}$ converges in distribution to $X$ if  for all bounded continuous functions $f \colon C^n[0,T] \rightarrow \mathbb{R}$ 
\[
\mathbb{E} f(\bar{X}_{\Delta t}) \rightarrow \mathbb{E} f(X),
\]
as $\Delta t \rightarrow 0$.  Stroock and Varadhan's result gives conditions on the original system of  SDEs and the numerical method under which $\bar{X}_{\Delta t}$ converges in distribution to $X$.

The system of SDEs (\ref{eqn:sde}) we consider  is determined by its coefficients $a$ and $\sigma_r$.  We define the matrix $b$ from $\sigma_r$ by
\begin{equation} \label{eqn:bdef}
b_{ij} (x) = \sum_{r} \sigma_{r,i}(x) \sigma_{r,j}(x),
\end{equation}
for $i,j=1,\ldots,n$.
Recall that the increment of the numerical method  (\ref{eqn:method}) starting from $x$ is denoted $\bar{\delta}(x,\Delta t)$.
For our numerical method we define corresponding coefficients $a_{\Delta t}$ and $b_{\Delta t}$ by 
\[
a_{\Delta t,i}(x) = \frac{1}{\Delta t} \mathbb{E} \bar{\delta}_i(x,\Delta t)  \mathbf{1}_{|\bar{\delta}(x, \Delta t)| \leq 1}, 
\]
and
\[
b_{\Delta t, ij}(x) = \frac{1}{\Delta t} \mathbb{E} \bar{\delta}_i(x,\Delta t) \bar{\delta}_j(x,\Delta t)  \mathbf{1}_{|\bar{\delta}(x, \Delta t)| \leq 1}.
\]
Finally, we define $\Gamma^{\epsilon}_{\Delta t}$ by
\[
\Gamma^{\epsilon}_{\Delta t}(x) = \frac{1}{\Delta t} \mathbb{P}( |\bar{\delta}(x, \Delta t)| \geq \epsilon ).
\]

In the following let $\| \cdot \|$ denote any norm on the space of $n \times n$ matrices.

\begin{thm}\label{thm:path}
Suppose that \\
(a) the coefficients $a$ and $\sigma_r$ 
of the SDEs (\ref{eqn:sde}) are locally Lipschitz continuous; \\
(b) 
there is a constant $C$ such that
for all $x \in \mathbb{R}^n$
\[
x^T a(x)  \leq C ( 1+ |x|^2) \ \ \  \mbox{and} \ \ \
\| b(x) \| \leq C (1+|x|^2);
\] 
(c) for all $R>0$
\[
\lim_{\Delta t \rightarrow 0} \sup_{|x| \leq R} | a_{\Delta t}(x) - a(x) | = 0,  \ \ \ \ \ \ 
\lim_{\Delta t \rightarrow 0} \sup_{|x| \leq R} \| b_{\Delta t}(x) - b(x) \| = 0, 
\]
and
\[
\lim_{\Delta t \rightarrow 0} \sup_{|x| \leq R} \Gamma^{\epsilon}_{\Delta t}(x) = 0,
\]
for all $\epsilon>0$.
\\
Then $\bar{X}_{\Delta t}$ converges in distribution to $X$ in $C^n[0,T]$.
\end{thm}

\begin{proof}
This result is Theorem 11.2.3 of \cite{SandV}, using Theorems 5.3.1 and 5.3.2 of  \cite{durrett} to obtain the well-posedness of the martingale problem.
\end{proof}

We remark that condition (b) of Theorem~\ref{thm:path}  is stronger than necessary.  The result holds with (b) replaced by the weaker condition that the martingale problem with coefficients $a$ and $b$ has a unique solution.

To give a feeling for the power of this result, here are some examples of functions $f$ for which it applies.  Firstly, we can recover the simpler pointwise results (without rates) if we let $f \colon C^n[0,T] \rightarrow  \mathbb{R}$
be defined by
\(
f(X) := g(X(t)),
\)
for some time $t \in [0,T]$.
If $g \colon \mathbb{R} \rightarrow \mathbb{R}$ is continuous and bounded then $f$ is continuous and bounded and the previous theorem tells us that $\mathbb{E}g(\bar{X}_{\Delta t}(t)) \rightarrow \mathbb{E} g(X(t))$.   More generally, we can choose $f$ to depend on $X$ through a number of points $t_1, \ldots, t_k \in [0,T]$.  Suppose $g \colon \mathbb{R}^k \rightarrow \mathbb{R}$ is bounded and continuous.  Letting $f(X) := g(X(t_1), \ldots, X(t_k) )$ gives us
\[
\mathbb{E} g(\bar{X}_{\Delta t}(t_1), \ldots, \bar{X}_{\Delta t}(t_k))
\rightarrow 
\mathbb{E} g(X(t_1), \ldots, X(t_k)),
\]
as $\Delta t \rightarrow 0$.  Even more generally, we can look at functions that do not depend on any finite number of times $t_i$.  For example, let $g \colon \mathbb{R} \rightarrow \mathbb{R}$ be bounded and continuous and define
\(
f(X) := g( \max_{t \in [0,T]} X(t)).
\)
The previous theorem tells us that  $\mathbb{E}f(\bar{X}_{\Delta t}) \rightarrow \mathbb{E}f(X)$ as $\Delta t \rightarrow 0$.

These example all rely on $f$ being continuous and bounded.  However, these assumptions can be weakened considerably in some cases.  One result of this type is that if $f$ is bounded and measurable and the probability of $X$ falling in the set of discontinuities of $f$ is zero, then it still holds that $\mathbb{E}f(\bar{X}_{\Delta t})  \rightarrow 
\mathbb{E}f(X)$ \cite[p.\ 21]{billingsley}.   As an example, let $f(X) = \mathbf{1}_{X(t_1)\in A_1} \mathbf{1}_{X(t_2) \in A_2}$.  Then if we can show $\mathbb{P}( X(t_1) \in \partial A_1, X(t_2) \in \partial A_2)=0$, it follows that 
\[
\mathbb{P}( \bar{X}_{\Delta t}(t_1) \in A_1,\bar{X}_{\Delta t}(t_2) \in A_2 ) 
\rightarrow 
\mathbb{P}( X(t_1) \in A_1,X(t_2) \in A_2 ) 
\]
as $\Delta t \rightarrow 0$.
Similarly, the condition on the boundedness of $f$ can be relaxed given some a priori knowledge on  the distributions of $X_n$ and  $X$ \cite[p.\ 31]{billingsley}. 

We now put the conditions of Theorem~\ref{thm:path} in terms more familiar in the numerical analysis of SDEs.

\begin{thm} \label{thm:mysde}
Suppose the coefficients $a$ and $\sigma$ of (\ref{eqn:sde}) are locally Lipschitz continuous and condition (b) of Theorem~\ref{thm:path} is satisfied.  Suppose that
the following limits hold as $\Delta t \rightarrow 0$ uniformly on bounded subsets of $\mathbb{R}^n$:
\begin{equation} \label{eqn:mycond1}
\frac{1}{\Delta t} \left| \mathbb{E} \delta_i - \mathbb{E} \bar{\delta}_i \right|  \rightarrow 0,
\end{equation}
\begin{equation} \label{eqn:mycond2}
\frac{1}{\Delta t} \left| \mathbb{E} \delta_i \delta_j - \mathbb{E} \bar{\delta}_i \bar{\delta}_j \right| \rightarrow 0,
\end{equation}
and 
\begin{equation} \label{eqn:mycond3}
\frac{1}{\Delta t}\mathbb{E} \left| \bar{\delta}_i \bar{\delta}_j \bar{\delta}_k \right| \rightarrow 0,
\end{equation}
for all $i,j,k= 1,\ldots,n$.
Then $\bar{X}_{\Delta t}$ converges in distribution to $X$ in $C^n[0,T]$.
\end{thm}
\begin{proof}
We prove this result by showing that  (\ref{eqn:mycond1}),  (\ref{eqn:mycond2}),  and (\ref{eqn:mycond3}) imply condition (c) of Theorem~\ref{thm:path}.  In the following we suppress the arguments of $\delta(x,\Delta t)$ and  $\bar{\delta}(x,\Delta t)$.   Fix a bounded set in $\mathbb{R}^n$. 
 We start by proving that (\ref{eqn:mycond3}) implies 
$\Gamma^{\epsilon}_{\Delta t} (x) \rightarrow 0$ uniformly.  Using Chebyshev's inequality:
\[
\Gamma^{\epsilon}_{\Delta t} (x)  =  \frac{1}{\Delta t} \mathbb{P} (| \bar{\delta}(x,\Delta t)| \geq \epsilon)  \leq   \sum_{i=1}^n\frac{1}{\Delta t} \mathbb{P}( | \bar{\delta}_i | \geq \epsilon/n ) 
 =  \sum_{i=1}^n \frac{1}{\Delta t} \mathbb{E} |\bar{\delta}_i| ^3  (n/\epsilon)^3,
\]
which goes to zero as $\Delta t \rightarrow 0$.

Next we prove that  $\| b_{\Delta t} (x)- b(x) \|$ goes to zero uniformly as $\Delta t \rightarrow 0$.
  The argument for $|a_{\Delta t}(x) -a(x) |$ is analogous so we omit it.
  Note that it suffices to show that $|b_{\Delta t, i, j} - b_{i,j} | \rightarrow 0$ uniformly on the chosen set for all $i,j$.
For each $i$ and $j$
\begin{align*}
| b_{\Delta t,ij} - b_{ij} | & =  
| \frac{1}{\Delta t}  \mathbb{E} \bar{\delta}_i \bar{\delta}_j \mathbf{1}_{|\bar{\delta}| \leq 1}
- b_{ij} | \\
& \leq    \frac{1}{\Delta t} | \mathbb{E} \bar{\delta}_i \bar{\delta}_j \mathbf{1}_{|\bar{\delta}| \leq 1}
- \mathbb{E} \bar{\delta}_i \bar{\delta}_j | + \frac{1}{\Delta t} | \mathbb{E} \bar{\delta}_i \bar{\delta}_j - \mathbb{E} \delta_i \delta_j | + | \frac{1}{\Delta t}\mathbb{E} \delta_i \delta_j - b_{ij} |
\end{align*}
The second term goes to zero from (\ref{eqn:mycond2}).  The third term goes to zero by general properties of SDEs.  The first term is equal to
\begin{align*}
\frac{1}{\Delta t} | \mathbb{E} \bar{\delta}_i \bar{\delta}_j \mathbf{1}_{|\bar{\delta}| >1} |
& \leq  \sum_{k=1}^n \frac{1}{\Delta t}  \mathbb{E}  |\bar{\delta}_i\bar{\delta}_j |
\mathbf{1}_{|\bar{\delta}_k| > 1/n} 
\leq  n \sum_{k=1}^n \frac{1}{\Delta t} \mathbb{E} |\bar{\delta}_i\bar{\delta}_j \bar{\delta}_k |,
\end{align*}
which goes to zero by (\ref{eqn:mycond3}).
\end{proof}

In order to show a strong re-embedding type of result, we use the following theorem, sometimes called Skorohod's Theorem.

\begin{thm} (See \cite[p.\ 70]{billingsley}.)
Let $S$ be a separable metric space with metric $d$.  Suppose that $X_n$, $n \geq 1$ and $X$ are random variables taking values in $S$, and that $X_n$ converges in distribution to $X$.  Then 
there are random variables $Y_n$, $n \geq 1$ and $Y$ all defined on the same probability space such that the distribution of $Y_n$ is the same as $X_n$ for all $n$, the distribution of $Y$ is the same as 
$X$, and  $Y_n$ converges to $Y$ almost surely.
\end{thm} 

Applying this theorem to $\bar{X}_{\Delta t}$, and $X$ gives the following result.

\begin{thm}  Let either the conditions of Theorem~\ref{thm:path} or the conditions of 
Theorem~\ref{thm:mysde} hold.
 Let $\Delta t_n$ be a sequence of positive step-sizes converging to $0$.  Then there are random elements $Y_n$ and $Y$ of $C^n[0,T]$ such that $Y_n$ has the same distribution as $\bar{X}_{\Delta t_n}$, $Y$ has the same distribution as $X$, and $Y_n \rightarrow Y$ in $C^n[0,T]$ almost surely.
  Thus, almost surely,
 \[
 \lim_{n \rightarrow \infty} \sup_{t \in [0,T]} | Y_n(t) - Y(t) | =0.
 \]
\end{thm}   \hfill $\square$

\vspace{1cm}

\noindent {\bf Acknowledgments.}   BC was funded by a Canadian National Science and Engineering Research Council (NSERC) postdoctoral fellowship.  YS was supported by 
an Institut des Sciences Math\'{e}matiques undergraduate research fellowship.  PFT was funded by an NSERC Discovery Grant.  The authors would like to thank the reviewers for comments that greatly improved the quality of the manuscript.

\bibliography{../../../paulsbib}{}
\bibliographystyle{abbrv}

\end{document}